\documentclass[a4paper, 12pt, reqno]{amsart}
\usepackage[textwidth=16cm, hmarginratio=1:1]{geometry}
\usepackage[T1]{fontenc}
\usepackage{amsmath,amssymb,amsthm,mathtools}
\usepackage{array,longtable,enumitem,microtype}
\usepackage{tikz}
\usepackage{mathrsfs}
\usepackage{float}
\usetikzlibrary{arrows.meta,positioning}
\usepackage[colorlinks]{hyperref} %backref=page
\allowdisplaybreaks

\hypersetup{
bookmarksnumbered,
pdfstartview={FitH},
breaklinks=true,
linkcolor=blue,
urlcolor=blue,
citecolor=blue,
bookmarksdepth=2
}

\newtheorem{theorem}{Theorem}[section]
\newtheorem{proposition}[theorem]{Proposition}
\newtheorem{lemma}[theorem]{Lemma}
\newtheorem{corollary}[theorem]{Corollary}
\theoremstyle{definition}
\newtheorem{definition}[theorem]{Definition}

\newcommand{\IAM}{\operatorname{IAM}}

\newcommand{\eps}{\varepsilon}

\newcommand{\Bool}{\mathcal B}

\newcommand{\IB}{\operatorname{IB}}

\usetikzlibrary{arrows.meta,decorations.markings}

\definecolor{strandA}{RGB}{35,92,170}
\definecolor{strandB}{RGB}{191,91,30}
\tikzset{
  edge/.style={draw=black!35,thin},
  neighbourhood/.style={draw=black!45,dashed,thin},
  motion/.style={thick,postaction={decorate},
    decoration={markings,mark=at position 0.55 with
      {\arrow{Latex[length=2.1mm,width=1.4mm]}}}},
  a-motion/.style={motion,draw=strandA},
  b-motion/.style={motion,draw=strandB},
  dot/.style={circle,inner sep=1.65pt},
  a-dot/.style={dot,fill=strandA},
  b-dot/.style={dot,fill=strandB},
  c-dot/.style={dot,fill=black}
}

\newcommand{\TwistRow}[4]{%
  \begin{scope}[yshift=#1]
    \node at (0,1.12) {#4};
    \draw[edge] (-1,0)--(1,0)--(0,-1.6)--cycle;
    \draw[neighbourhood] (0,0) ellipse[x radius=1.48,y radius=.72];
    \node[font=\small] at (-1.15,.64) {};
    \ifnum#2=1
      \draw[a-motion] (-1,0)
        arc[start angle=180,end angle=0,x radius=1,y radius=.46];
      \draw[b-motion] (1,0)
        arc[start angle=0,end angle=-180,x radius=1,y radius=.46];
    \else
      \draw[a-motion] (-1,0)
        arc[start angle=180,end angle=360,x radius=1,y radius=.46];
      \draw[b-motion] (1,0)
        arc[start angle=0,end angle=180,x radius=1,y radius=.46];
    \fi
    \node[a-dot,label=left:$A$] at (-1,0) {};
    \node[b-dot,label=right:$B$] at (1,0) {};
    \node[c-dot,label=below:$C$] at (0,-1.6) {};
    \node[anchor=west,font=\small] at (.17,-1.6) {fixed};
    \draw[-{Latex[length=2.7mm]},semithick]
      (1.9,-.18)--node[above=3pt] {$#3$} (3.25,-.18);
    \begin{scope}[xshift=5.2cm]
      \node at (0,1.12) {After the half twist};
      \draw[edge] (-1,0)--(1,0);
      \draw[neighbourhood] (0,0)
        ellipse[x radius=1.48,y radius=.72];
      \node[b-dot,label=left:$B$] at (-1,0) {};
      \node[a-dot,label=right:$A$] at (1,0) {};
      \node[c-dot,label=below:$C$] at (0,-1.6) {};
      \node[anchor=west,font=\small] at (.17,-1.6) {fixed};
    \end{scope}
  \end{scope}%
}

\title{Presentations of inverse Artin monoids of type $D_n$ via cluster mutations}
\author{Bing Duan$^\dag$}
\address[Bing Duan]{School of Mathematics and Statistics, Lanzhou University, Lanzhou, 730000 P. R. China}
\email{duanbing@lzu.edu.cn}
\thanks{$\dag$ Corresponding author}
\date{}

\begin{document}

\maketitle

\begin{abstract}
Let $D_n^\varepsilon$ denote the Dynkin diagram of type $D_{n+1}$ with the endpoint $\varepsilon$ of its long branch frozen, and let $\mathcal{M}(D_n^\varepsilon)$ be the mutation class of any orientation of this diagram under mutations at mutable vertices. For each $\Gamma\in\mathcal{M}(D_n^\varepsilon)$, we define an inverse monoid $M_\Gamma$ by generators and relations determined by $\Gamma$. We prove that $M_\Gamma\cong M_{\mu_k(\Gamma)}$ for every mutable vertex $k$. Consequently, every quiver in this mutation class yields a presentation of the inverse Artin monoid of type $D_n$. Our argument establishes mutation invariance for the presentations of inverse Artin monoids of both types $A_n^\varepsilon$ and $D_n^\varepsilon$.

\hspace{0.15cm}


\end{abstract}

\section{Introduction}\label{sec:introduction}

Braid groups, introduced by Artin \cite{Artin1925,Artin1947}, are fundamental objects in algebra and topology. Artin groups generalize braid groups and, in finite type, have the corresponding finite reflection groups as quotients. The inverse braid monoid introduced by Easdown and Lavers \cite{EL2004} extends the braid group by allowing strands to be deleted. Everitt and Fountain \cite{EF10,EF13} introduced the notion of Boolean reflection monoids, which extend finite reflection groups by incorporating partial symmetries. Building on this framework, Duan~\cite{Duan2025} introduced inverse Artin monoids associated with Boolean reflection monoids and constructed their partial braid realizations.

Cluster algebras were introduced by Fomin and Zelevinsky in \cite{FZ2002}. Surprisingly, the classification of cluster algebras of finite type, that is, those with finitely many cluster variables, coincides with the Cartan--Killing classification of complex simple Lie algebras \cite{FZ2003}. Barot and Marsh \cite{BM2015} obtained presentations of finite crystallographic reflection groups from arbitrary seeds of the corresponding cluster algebras. Their presentations supplement the usual relations between pairs of generators with relations associated with chordless cycles. Analogous presentations of Artin groups were obtained by Haley, Hemminger, Landesman and Peck~\cite{HHLP2017} and, in the simply laced setting, by Grant and Marsh~\cite{GM2017}. The associated mutation isomorphisms are defined using conjugations determined by the arrows incident to the mutated vertex. These results show how quiver mutation can be used to construct presentations of reflection groups, Artin groups, and related algebraic structures.

More recently, Fedele and Marsh~\cite{FM2026} obtained compatible presentations of the complex reflection groups $G(d,d,n)$ and their braid groups for $d\geq2$. These presentations are indexed by tagged triangulations of a disk with $n$ boundary marked points and a single interior cone point of degree $d$. Their construction generalizes the type $D_n$ presentations of Grant and Marsh, which are recovered in the case $d=2$.

The successful interaction between cluster algebras and presentations of groups motivates our study of cluster algebras and Boolean reflection monoids \cite{EF10,EF13}, as well as inverse Artin monoids \cite{Duan2025}, which are inverse semigroups whose properties are very close to those of reflection groups and Artin groups. Duan, Li and Luo \cite{DLL2020} obtained presentations of Boolean reflection monoids compatible with mutations of quivers carrying a frozen vertex. Unlike the case of groups, when dealing with presentations of inverse semigroups, we need an additional idempotent generator, which suggests the use of quivers with frozen vertices. These presentations of Boolean reflection monoids provide the combinatorial basis for the cycle and path relations used in this paper.

Let $D_n^\varepsilon$ be the Dynkin diagram of type $D_{n+1}$ with the terminal vertex at the end of its long branch frozen, and let $\mathcal{M}(D_n^\varepsilon)$ be the mutation class of an orientation of this diagram, with mutations allowed only at mutable vertices. For each $\Gamma\in\mathcal{M}(D_n^\varepsilon)$, we define a monoid $M_\Gamma$ by generators and relations determined by $\Gamma$ (see Definition \ref{def:B-inverse-monoid}), where each mutable vertex corresponds to a unit generator, and the frozen vertex corresponds to an idempotent. Our main result, Theorem \ref{thm:B-mutation-invariance}, states that $M_\Gamma\cong M_{\mu_k(\Gamma)}$ for every mutable vertex $k$. For an initial Dynkin orientation, our defining relations coincide with the standard presentation of $\operatorname{IAM}(D_n^\eps)$ in \cite{Duan2025}. Hence $M_\Gamma\cong\operatorname{IAM}(D_n^\eps)$ for every $\Gamma\in\mathcal M(D_n^\eps)$. Our proof combines the mutation results for Boolean reflection monoids in \cite[Definition 4.5 and Theorems 4.7 and 4.9]{DLL2020} with the idempotent-separating property of the canonical epimorphism from $M_\Gamma$ to the Boolean reflection monoid of type $D_n$, see Proposition \ref{idemsep}. Our approach to proving mutation invariance for type $D_n^\varepsilon$ differs from those in \cite{Duan2025,DL25}. In fact, our argument establishes mutation invariance for the presentations of inverse Artin monoids of both types $A_n^\varepsilon$ and $D_n^\varepsilon$.

The paper is organized as follows. In Section~\ref{sec:preliminaries}, we review some basic facts on ice quivers and their mutations, together with the relevant inverse Artin monoids. In Section~\ref{monoids and diagrams}, we first define a presentation associated with each quiver in $\mathcal{M}(D_n^\varepsilon)$, and then state our main result, Theorem~\ref{thm:B-mutation-invariance}. Finally, in Section~\ref{proof of main theorem}, we prove Theorem~\ref{thm:B-mutation-invariance}.

\section{Preliminaries}\label{sec:preliminaries}

\subsection{Ice quivers and their mutations}\label{Ice quivers and their mutations}

Let $I$ be a finite set with $|I|=n$, and let $\varepsilon\notin I$. An ice quiver $Q$ is a finite quiver with vertex set $I\cup\{\varepsilon\}$, with no loops or oriented $2$-cycles. The vertices in $I$ are mutable, and $\varepsilon$ is frozen. For a mutable vertex $k\in I$, the mutation $\mu_k(Q)$ is defined by the following rule of Fomin and Zelevinsky in \cite{FZ2002}:
\begin{enumerate}
  \item[(1)] For each oriented path $i\longrightarrow k\longrightarrow j$ in $Q$, add an arrow $i\longrightarrow j$.

  \item[(2)] Reverse every arrow incident with $k$.

  \item[(3)] Remove a maximal collection of oriented $2$-cycles.
\end{enumerate}
The vertices $i$ and $j$ in (1) may include $\varepsilon$. Thus, arrows incident with the frozen vertex are subject to the same mutation rule, although mutation at $\varepsilon$ is not allowed. Mutation is involutive: $\mu_k^2(Q)=Q$. %%%%For simplicity, we write $\mathcal{M}(Q)$ for the class of quivers obtained from $Q$ by finite sequences of mutations at mutable vertices, considered up to isomorphisms preserving the distinction between mutable and frozen vertices. 

Let $T$ be a finite tree and fix a vertex $v$. As shown in the proof of Proposition 9.2 of Fomin and Zelevinsky in \cite{FZ2003}, any two orientations of $T$ can be connected by a sequence of mutations at sources or sinks avoiding $v$. The argument proceeds by removing a leaf $\ell\ne v$ and applying induction. The resulting sequence lifts to $T$ by inserting mutations at $\ell$ whenever necessary to ensure that its neighbour remains a source or sink when mutated. Finally, a mutation at $\ell$, if needed, gives the remaining edge its desired orientation. Consequently, when $v$ is frozen and all other vertices are mutable, the mutation class is independent of the initial orientation of $T$.

Let $\mathcal{M}(A_n)$ denote the mutation class of an orientation of the Dynkin diagram of type $A_n$, with all vertices mutable. The diagrams $A_n^\varepsilon$ and $D_n^\varepsilon$ are obtained by freezing an endpoint of $A_{n+1}$ and the endpoint of the long branch of $D_{n+1}$, respectively; see Figure~\ref{the initial diagrams}. Thus, in both cases, $n$ is the number of mutable vertices. We denote the corresponding mutation classes of ice quivers by $\mathcal{M}(A_n^\varepsilon)$ and $\mathcal{M}(D_n^\varepsilon)$, with mutations allowed only at mutable vertices. The preceding argument shows that the mutation classes $\mathcal{M}(A_n^\varepsilon)$ and $\mathcal{M}(D_n^\varepsilon)$ are independent of the initial orientation.

\tikzset{
  iam/vertex/.style={
    circle,
    draw,
    fill=white,
    inner sep=0pt,
    minimum size=4pt
  },
  iam/frozen/.style={
    iam/vertex,
    fill=black
  },
  iam/arrow/.style={
    ->,
    >=stealth,
    shorten >=1pt,
    shorten <=1pt
  },
  iam/edge/.style={
    shorten >=1pt,
    shorten <=1pt
  }
}

\begin{figure}[htbp]
\centering
\begin{tikzpicture}[x=1.3cm,y=1cm,font=\small]
  \node[left] at (-0.45,0) {$A_n^\varepsilon$};

  \node[iam/vertex,label=below:$1$] (a1) at (0,0) {};
  \node[iam/vertex,label=below:$2$] (a2) at (1,0) {};
  \node[iam/vertex,label=below:$n$] (an) at (3,0) {};
  \node[iam/frozen,label=below:$\varepsilon$] (ae) at (4,0) {};

  \draw[iam/edge] (a1)--(a2);
  \draw[iam/edge,densely dashed] (a2)--(an);
  \draw[iam/edge] (an)--(ae);

  \node[left] at (-0.45,-2.1) {$D_n^\varepsilon$};

  \node[iam/vertex,label=left:$1$] (d0) at (0,-2.7) {};
  \node[iam/vertex,label=left:$2$] (d1) at (0,-1.5) {};
  \node[iam/vertex,label=below:$3$] (d2) at (1,-2.1) {};
  \node[iam/vertex,label=below:$n$] (dn) at (3,-2.1) {};
  \node[iam/frozen,label=below:$\varepsilon$] (de) at (4,-2.1) {};

  \draw[iam/edge] (d0)--(d2);
  \draw[iam/edge] (d1)--(d2);
  \draw[iam/edge,densely dashed] (d2)--(dn);
  \draw[iam/edge] (dn)--(de);
\end{tikzpicture}
\caption{The diagrams $A_n^\varepsilon$ and $D_n^\varepsilon$ ($n\geq4$). The filled vertex is frozen; the remaining vertices are mutable.} \label{the initial diagrams}
\end{figure}
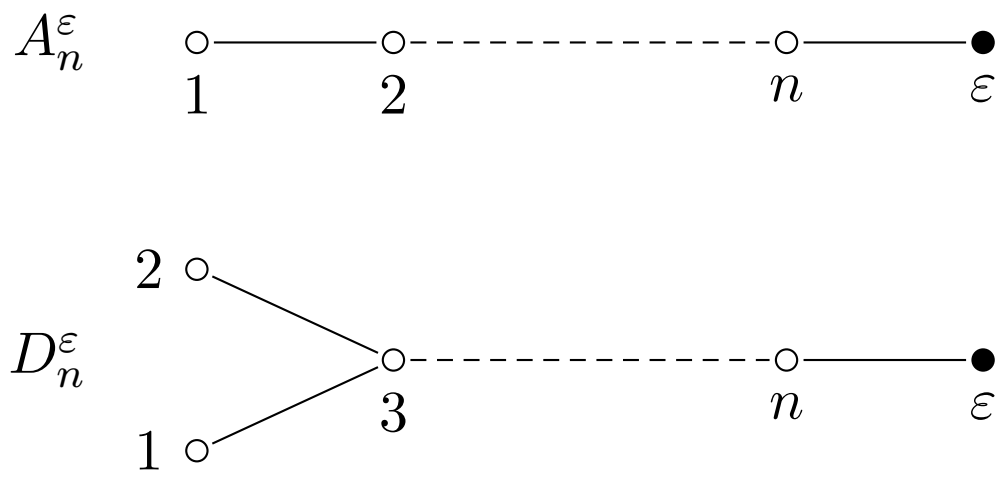

The mutation classes $\mathcal{M}(A_n^\varepsilon)$ and $\mathcal{M}(D_n^\varepsilon)$ are finite. Indeed, unfreezing the distinguished vertex in each initial quiver yields a quiver of Dynkin type $A_{n+1}$ or $D_{n+1}$, respectively, whose mutation class is finite by the finite type classification of Fomin and Zelevinsky \cite[Theorem 1.4]{FZ2003}. The assertion follows since freezing a vertex preserves finiteness of the mutation class.

We next recall the descriptions of the mutation classes $\mathcal{M}(A_n^\varepsilon)$ and $\mathcal{M}(D_n^\varepsilon)$ following \cite[Section~4.2]{DLL2020} and \cite[Section~2]{Vatne2010}. These descriptions can also be obtained from the polygon models presented in \cite[Chapter~3]{Schiffler2014}. The quivers in the descriptions below are connected and have at most one arrow between any two vertices. A chordless cycle is a cycle of length at least three in the underlying undirected graph with no edge joining nonconsecutive vertices. It is oriented if its arrows all follow the same cyclic direction. The degree of a vertex is the number of its neighbours.

Let $Q$ be a connected ice quiver with $n$ mutable vertices, a unique frozen vertex $\varepsilon$, and at most one arrow between any two vertices. Then $Q\in\mathcal{M}(A_n^\varepsilon)$ if and only if the following conditions hold:
\begin{enumerate}
  \item[(1)] Every chordless cycle is an oriented $3$-cycle.

  \item[(2)] Every mutable vertex has degree at most four.

  \item[(3)] If a mutable vertex has degree four, its incident arrows can be partitioned into two pairs, each contained in an oriented $3$-cycle.

  \item[(4)] If a mutable vertex has degree three, two of its incident arrows lie in a common oriented $3$-cycle, and the remaining arrow lies in no oriented $3$-cycle.

  \item[(5)] The frozen vertex $\varepsilon$ has degree at most two; if its degree is two, it lies in an oriented $3$-cycle.
\end{enumerate}

A vertex of a quiver of type $A$ is called a connecting vertex if it has degree at most two and, in the case of degree two, lies in an oriented $3$-cycle. In particular, the frozen vertex $\varepsilon$ is a connecting vertex in every quiver in $\mathcal{M}(A_n^\varepsilon)$. We adopt the conventions that $\mathcal{M}(A_1)$ consists of the quiver with a single mutable vertex and $\mathcal{M}(A_0^\varepsilon)$ consists of the quiver with the single frozen vertex $\varepsilon$. For a set $S$ of vertices of $Q$, we write $Q\setminus S$ for the full subquiver obtained by deleting $S$ and all incident arrows. When an arrow is deleted, both endpoints are retained.

\begin{figure}[htbp]
\centering
\begin{tikzpicture}[font=\small]
  \draw (1.5,0) ellipse (1.8cm and 0.85cm);

  \node[iam/vertex,label=above:$a$] (a) at (-1.3,1) {};
  \node[iam/vertex,label=below:$b$] (b) at (-1.3,-1) {};
  \node[iam/vertex,label=below:$c$] (c) at (0,0) {};
  \node[iam/frozen,label=right:$\varepsilon$] (e) at (2.4,0) {};

  \draw[iam/edge] (a)--(c);
  \draw[iam/edge] (b)--(c);
  \draw[iam/edge,densely dashed] (c)--(e);

  \node at (1.5,0.4) {$Q'$};
\end{tikzpicture}
\caption{Type I. The two edges outside $Q'$ may have either orientation. The dashed line represents the unspecified part of $Q'$.} \label{fig:type-I}
\end{figure}
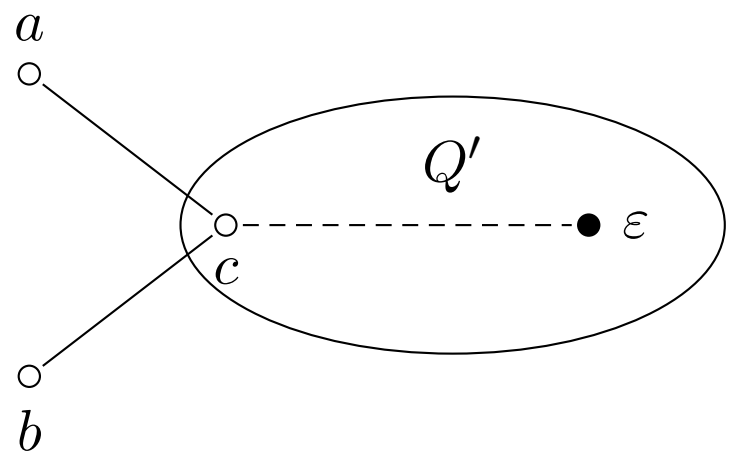

In the classification below, all conditions on connecting vertices are understood relative to the indicated component. A connected ice quiver $Q$ with $n$ mutable vertices and a unique frozen vertex $\varepsilon$ belongs to $\mathcal{M}(D_n^\varepsilon)$ if and only if it has one of the following four forms:

\medskip
\noindent\textbf{Type I.}
There are distinct mutable vertices $a$, $b$ and $c$ such that $a$ and $b$ have degree one and their common neighbour is $c$. The arrows between $c$ and $a,b$ may be oriented independently. Moreover,
\[
Q'=Q\setminus\{a,b\}\in\mathcal{M}(A_{n-2}^\varepsilon),
\]
and $c$ is a connecting vertex of $Q'$; see Figure~\ref{fig:type-I}. In particular, the frozen vertex belongs to $Q'$ and is distinct from $c$.

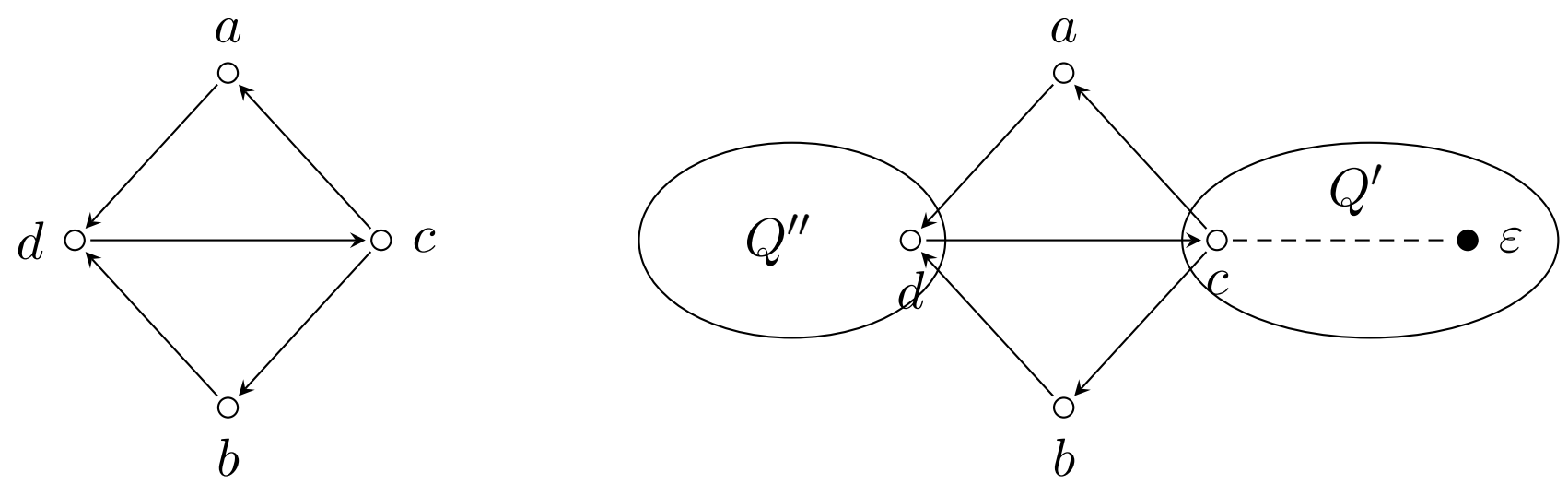
\begin{figure}[htbp]
\centering
\begin{tikzpicture}[font=\small]
  \node[iam/vertex,label=left:$d$] (d) at (0,0) {};
  \node[iam/vertex,label=right:$c$] (c) at (2.2,0) {};
  \node[iam/vertex,label=above:$a$] (a) at (1.1,1.2) {};
  \node[iam/vertex,label=below:$b$] (b) at (1.1,-1.2) {};

  \draw[iam/arrow] (d)--(c);
  \draw[iam/arrow] (c)--(a);
  \draw[iam/arrow] (a)--(d);
  \draw[iam/arrow] (c)--(b);
  \draw[iam/arrow] (b)--(d);

  \begin{scope}[xshift=6cm]
    \draw (-0.85,0) ellipse (1.1cm and 0.7cm);
    \draw (3.3,0) ellipse (1.35cm and 0.7cm);

    \node[iam/vertex,label=below:$d$] (dd) at (0,0) {};
    \node[iam/vertex,label=below:$c$] (cc) at (2.2,0) {};
    \node[iam/vertex,label=above:$a$] (aa) at (1.1,1.2) {};
    \node[iam/vertex,label=below:$b$] (bb) at (1.1,-1.2) {};
    \node[iam/frozen,label=right:$\varepsilon$] (ee) at (4,0) {};

    \draw[iam/arrow] (dd)--(cc);
    \draw[iam/arrow] (cc)--(aa);
    \draw[iam/arrow] (aa)--(dd);
    \draw[iam/arrow] (cc)--(bb);
    \draw[iam/arrow] (bb)--(dd);
    \draw[iam/edge,densely dashed] (cc)--(ee);

    \node at (-0.95,0) {$Q''$};
    \node at (3.2,0.36) {$Q'$};
  \end{scope}
\end{tikzpicture}
\caption{Type II: the four-vertex subquiver (left) and its attached components (right). The opposite orientation is also allowed.} \label{fig:type-II}
\end{figure}

\medskip
\noindent\textbf{Type II.}
The full subquiver on four distinct vertices $a,b,c,d$ consists of the two oriented $3$-cycles
\[
d\longrightarrow c\longrightarrow a\longrightarrow d, \qquad d\longrightarrow c\longrightarrow b\longrightarrow d,
\]
sharing the arrow $d\longrightarrow c$, or of the quiver obtained by reversing all five arrows. The vertices $a$ and $b$ are mutable and have degree two in $Q$. After deleting $a,b$ and the shared arrow, the remaining quiver is
\[
\bigl(Q\setminus\{a,b\}\bigr) \setminus\{d\longrightarrow c\}=Q'\sqcup Q'',
\]
with the reversed arrow deleted in the opposite orientation. Here $Q'$ and $Q''$ are the two connected components, with
\[
Q'\in\mathcal{M}(A_k^\varepsilon), \qquad Q''\in\mathcal{M}(A_\ell), \qquad
k\geq0,\quad \ell\geq1,\quad k+\ell+2=n.
\]
The vertex $c$ is a connecting vertex of $Q'$, and $d$ is a connecting vertex of $Q''$; see Figure~\ref{fig:type-II}. The frozen vertex lies in $Q'$. If $k=0$, then $Q'=\{\varepsilon\}$ and $c=\varepsilon$; otherwise, $c$ is mutable.

\begin{figure}[htbp]
\centering
\begin{tikzpicture}[font=\small]
  \node[iam/vertex,label=left:$d$] (d) at (0,0) {};
  \node[iam/vertex,label=right:$c$] (c) at (2.2,0) {};
  \node[iam/vertex,label=above:$a$] (a) at (1.1,1.2) {};
  \node[iam/vertex,label=below:$b$] (b) at (1.1,-1.2) {};

  \draw[iam/arrow] (c)--(a);
  \draw[iam/arrow] (a)--(d);
  \draw[iam/arrow] (d)--(b);
  \draw[iam/arrow] (b)--(c);

  \begin{scope}[xshift=6cm]
    \draw (-0.85,0) ellipse (1.1cm and 0.7cm);
    \draw (3.3,0) ellipse (1.35cm and 0.7cm);

    \node[iam/vertex,label=below:$d$] (dd) at (0,0) {};
    \node[iam/vertex,label=below:$c$] (cc) at (2.2,0) {};
    \node[iam/vertex,label=above:$a$] (aa) at (1.1,1.2) {};
    \node[iam/vertex,label=below:$b$] (bb) at (1.1,-1.2) {};
    \node[iam/frozen,label=right:$\varepsilon$] (ee) at (4,0) {};

    \draw[iam/arrow] (cc)--(aa);
    \draw[iam/arrow] (aa)--(dd);
    \draw[iam/arrow] (dd)--(bb);
    \draw[iam/arrow] (bb)--(cc);
    \draw[iam/edge,densely dashed] (cc)--(ee);

    \node at (-0.95,0) {$Q''$};
    \node at (3.2,0.36) {$Q'$};
  \end{scope}
\end{tikzpicture}
\caption{Type III: the oriented $4$-cycle (left) and its attached
components (right). The opposite orientation is also allowed.}
\label{fig:type-III}
\end{figure}
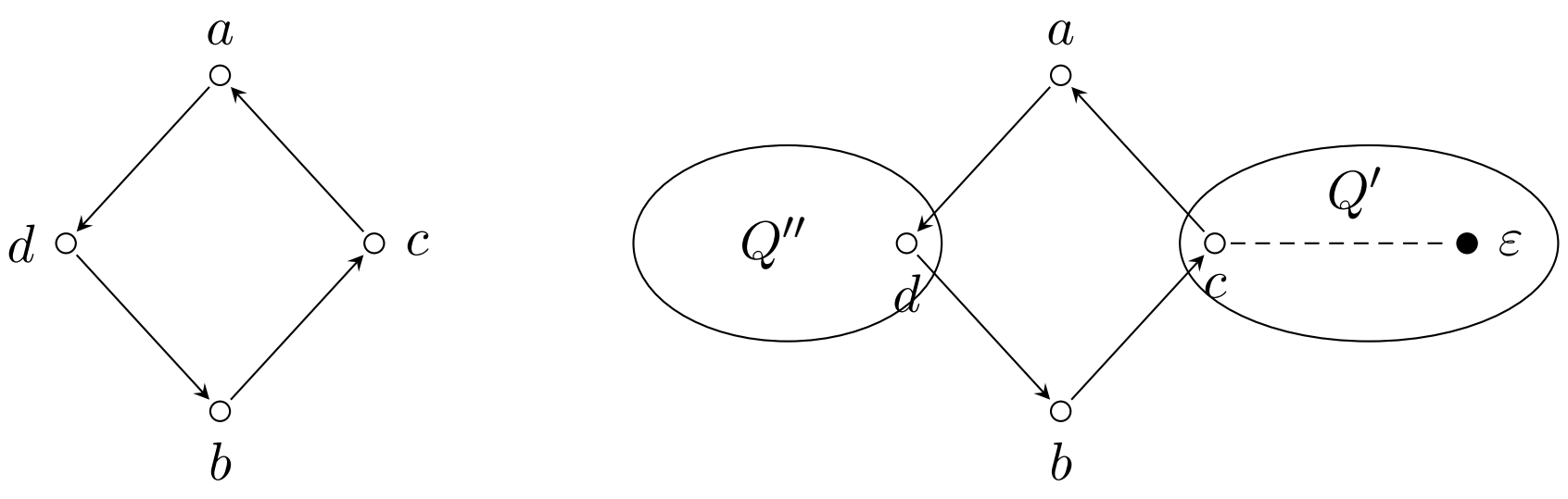

\medskip
\noindent\textbf{Type III.}
The full subquiver on four distinct vertices $a,b,c,d$ is the oriented $4$-cycle
$$
c\longrightarrow a\longrightarrow d\longrightarrow b
\longrightarrow c,
$$
or its reverse.
The vertices $a$ and $b$ are mutable and have degree two in $Q$.
Deleting them gives
$$
Q\setminus\{a,b\}=Q'\sqcup Q'',
$$
where $Q'$ and $Q''$ are the two connected components and
$$
Q'\in\mathcal{M}(A_k^\varepsilon),
\qquad
Q''\in\mathcal{M}(A_\ell),
\qquad
k\geq0,\quad \ell\geq1,\quad k+\ell+2=n.
$$
The vertices $c$ and $d$ are connecting vertices of $Q'$ and $Q''$, respectively; see Figure~\ref{fig:type-III}. As in Type II, $\varepsilon$ lies in $Q'$; $c=\varepsilon$ when $k=0$, and $c$ is mutable when $k>0$.

\medskip
\noindent\textbf{Type IV.}
There is a full subquiver consisting of an oriented cycle of length $d\geq3$, called the central cycle; all its vertices are mutable. For each arrow $\alpha:u\longrightarrow v$ of this cycle, there may be one additional vertex $c_\alpha$ outside the cycle such that the full subquiver on $u,v,c_\alpha$ is
\[
u\longrightarrow v\longrightarrow c_\alpha\longrightarrow u.
\]
We call this oriented $3$-cycle a spike and $c_\alpha$ its tip. Different spikes have distinct tips. Apart from the arrows of the central cycle and its spikes, no arrow is incident with a vertex of the central cycle. Deleting the vertices of the central cycle leaves a disjoint union
of connected components, one for each spike:
\[
Q\setminus\{\text{vertices of the central cycle}\}=Q'\sqcup Q_1\sqcup\cdots\sqcup Q_r.
\]
Exactly one component, denoted by $Q'$, contains $\varepsilon$, and
\[
Q'\in\mathcal{M}(A_k^\varepsilon), \qquad Q_j\in\mathcal{M}(A_{\ell_j}) \quad(1\leq j\leq r),
\]
where $k\geq0$, $r\geq0$, $\ell_j\geq1$, and $d+k+\sum_{j=1}^{r}\ell_j=n.$

\begin{figure}[htbp]
\centering
\begin{tikzpicture}[font=\small]
  \draw (0.85,3.5) ellipse (1.35cm and 1.05cm);
  \draw (-4,0.6) ellipse (1.35cm and 0.75cm);
  \draw (5.7,0.6) ellipse (1.35cm and 0.75cm);

  \node[iam/vertex,label=below right:$i_1$]
    (i1) at (-1.4,1.4) {};
  \node[iam/vertex,label=left:$i_2$]
    (i2) at (-1.4,-0.2) {};
  \node[iam/vertex,label=below:$i_3$]
    (i3) at (0,-1.3) {};
  \node[iam/vertex,label=below:$i_{d-2}$]
    (im2) at (1.7,-1.3) {};
  \node[iam/vertex,label=left:$i_{d-1}$]
    (im1) at (3.1,-0.2) {};
  \node[iam/vertex,label=below left:$i_d$]
    (id) at (3.1,1.4) {};

  \node[iam/vertex,label=below:$c'$]
    (c1) at (0.85,2.6) {};
  \node[iam/vertex,label=left:$c''$]
    (c2) at (-3,0.6) {};
  \node[iam/vertex,label=right:$c'''$]
    (c3) at (4.7,0.6) {};
  \node[iam/frozen,label=right:$\varepsilon$]
    (e) at (0.85,4) {};

  \draw[iam/arrow] (i1)--(i2);
  \draw[iam/arrow] (i2)--(i3);
  \draw[iam/arrow,densely dashed] (i3)--(im2);
  \draw[iam/arrow] (im2)--(im1);
  \draw[iam/arrow] (im1)--(id);
  \draw[iam/arrow] (id)--(i1);

  \draw[iam/arrow] (i1)--(c1);
  \draw[iam/arrow] (c1)--(id);

  \draw[iam/arrow] (i2)--(c2);
  \draw[iam/arrow] (c2)--(i1);

  \draw[iam/arrow] (id)--(c3);
  \draw[iam/arrow] (c3)--(im1);

  \draw[iam/edge,densely dashed] (c1)--(e);

  \node at (0.25,3.45) {$Q'$};
  \node at (-4.1,0.6) {$Q_1$};
  \node at (5.8,0.6) {$Q_2$};
\end{tikzpicture}
\caption{Type IV, illustrated with three spikes. The arrows specify the central cycle and the spikes;
the dashed portions and ellipses represent the remaining vertices and attached components.}
\label{fig:type-IV}
\end{figure}
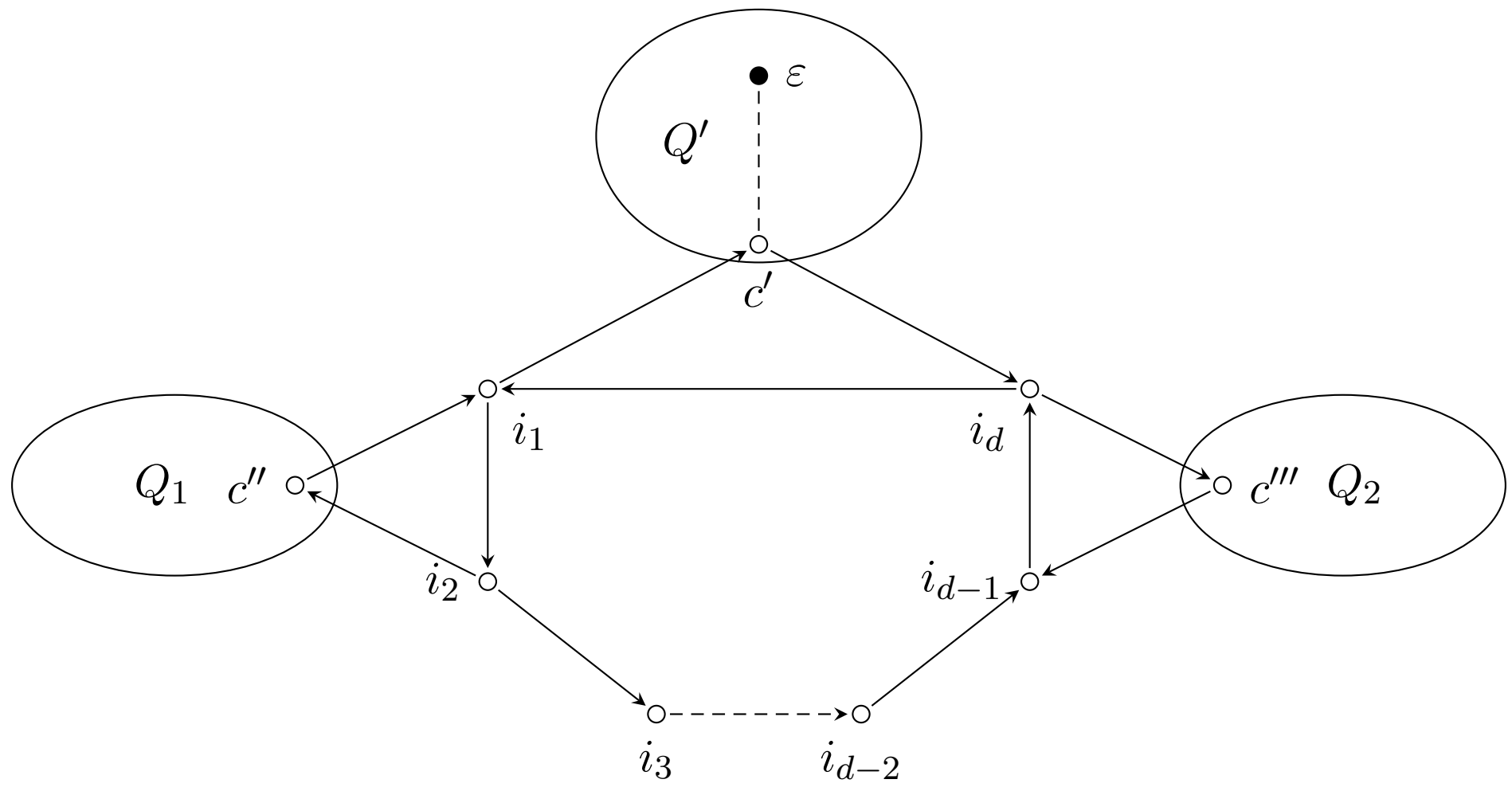

Each component contains precisely one spike tip, which is a connecting vertex of that component; see Figure~\ref{fig:type-IV}. In particular, there is at least one spike, since $\varepsilon$ is not on the central cycle. If $k=0$, then $Q'=\{\varepsilon\}$ and its spike tip is $\varepsilon$ itself. If $k>0$, that tip is mutable.

\subsection{Presentations of Artin groups}

Let $Q$ be a quiver with vertex set $I$ that is mutation equivalent to an orientation of a simply laced Dynkin diagram $\Delta$. An oriented chordless cycle is specified by an ordered tuple $C=(i_0,\ldots,i_{d-1})$ of distinct vertices such that
\[
i_0\longrightarrow i_1\longrightarrow\cdots
\longrightarrow i_{d-1}\longrightarrow i_0,
\qquad d\geq3,
\]
and the full subquiver on these vertices contains exactly the displayed arrows. For group elements $u,v$, let $[u,v]=uvu^{-1}v^{-1}$.

Following \cite[Definitions~4.4 and~4.5]{HHLP2017}, let $A_Q$ be the group generated by $x_i$, $i\in I$, subject to the following relations. For distinct vertices $i,j\in I$, impose
\[
\begin{aligned}
x_ix_j&=x_jx_i
&&\text{if $i$ and $j$ are not adjacent},\\
x_ix_jx_i&=x_jx_ix_j
&&\text{if $i$ and $j$ are adjacent}.
\end{aligned}
\]
For each oriented chordless cycle $C=(i_0,\ldots,i_{d-1})$, we impose
\[
[x_{i_a},p_a(C;x)]=1 \qquad(0\leq a<d),
\]
where all subscripts on the vertices of $C$ are read modulo $d$ and
\[
p_a(C;x)=x_{i_{a+1}}^{-1}\cdots x_{i_{a+d-2}}^{-1} x_{i_{a+d-1}} x_{i_{a+d-2}}\cdots x_{i_{a+1}}.
\]

\begin{theorem}[{\cite[Proposition~6.2 and Theorem~6.1]{HHLP2017}}]\label{thm:mutable}
If $Q'=\mu_k(Q)$, then
\[
\Phi_k:A_{Q'}\xrightarrow{\sim}A_Q,
\qquad x_i'\longmapsto
\begin{cases}x_kx_ix_k^{-1},&i\to k\text{ in }Q,\\x_i,&\text{otherwise},\end{cases}
\]
has inverse $x_i\mapsto (x_k')^{-1}x_i'x_k'$ when $i\to k$ in $Q$, and $x_i\mapsto x_i'$ otherwise. In particular, the group $A_Q$ is isomorphic to the Artin group of type $\Delta$. 
\end{theorem}

These conventions differ from the positive cycle convention of \cite[Definition 2.2]{GM2017} by inversion of the generators. Indeed, let $B_Q$ be the group generated by $y_i$, $i\in I$, with the same pair relations as $A_Q$, with $x_i$ replaced by $y_i$, and the following positive cycle relations.

For each oriented chordless cycle $C=(i_0,\ldots,i_{d-1})$, define
\[
W_a(C;y)=y_{i_a}y_{i_{a+1}}\cdots y_{i_{a+d-1}} y_{i_a}y_{i_{a+1}}\cdots y_{i_{a+d-3}}, \qquad 0\le a<d,
\]
where all subscripts on the vertices of $C$ are read modulo $d$. We impose
\[
 W_0(C;y)=W_1(C;y)=\cdots=W_{d-1}(C;y).
\]
These words are obtained from $W_0(C;y)$ by cyclically permuting the vertex labels along $C$.

Comparing Theorem~\ref{thm:mutable} with the inverse mutation isomorphism in \cite[Remark~2.11]{GM2017}, an induction along a mutation sequence from an orientation of $\Delta$ to $Q$ shows that the assignment
\[
\iota_Q:A_Q\longrightarrow B_Q,\qquad  x_i\longmapsto y_i^{-1},
\]
extends to an isomorphism.

%%The assignment
%%\begin{align*}
%%\iota_Q:A_Q & \longrightarrow B_Q, \quad x_i \longmapsto y_i^{-1}
%%\end{align*}
%%extends to an isomorphism. We prove this by induction along a mutation sequence from an orientation of $\Delta$ to $Q$. For the initial Dynkin orientation, there are no cycles, so both presentations consist only of the commutation and braid relations. Taking inverses of these relations shows that the assignments $x_i\mapsto y_i^{-1}$ and $y_i\mapsto x_i^{-1}$ define mutually inverse homomorphisms.

%%For the induction step, suppose that $\iota_Q$ is an isomorphism with the stated images of the generators, and let $Q'=\mu_k(Q)$. By \cite[Theorem~2.10 and Remark~2.11]{GM2017}, there is an isomorphism
%%$$
%%\begin{aligned}
%%H_k:B_{Q'}&\longrightarrow B_Q,\\
%%y_i'&\longmapsto
%%\begin{cases}
%%y_k^{-1}y_iy_k,
%%&\text{if $i\longrightarrow k$ in $Q$},\\
%%y_i,
%%&\text{otherwise}.
%%\end{cases}
%%\end{aligned}
%%$$
%%Together with the isomorphism $\Phi_k$ from Theorem~\ref{thm:mutable}, this defines an isomorphism
%%$$
%%\iota_{Q'}=H_k^{-1}\circ\iota_Q\circ\Phi_k:
%%A_{Q'}\longrightarrow B_{Q'}.
%%$$
%%If $i\longrightarrow k$ in $Q$, then
%%$$
%%\begin{aligned}
%%\iota_Q\bigl(\Phi_k(x_i')\bigr) =\iota_Q(x_kx_ix_k^{-1})=y_k^{-1}y_i^{-1}y_k=H_k\bigl((y_i')^{-1}\bigr).
%%\end{aligned}
%%$$
%%Otherwise, both sides are $y_i^{-1}$. Hence $\iota_{Q'}(x_i')=(y_i')^{-1}$ for every $i\in I$, which completes the induction.

\subsection{Inverse Artin monoid of type $D_n$}\label{sec:nverse Artin monoid of type Dn}

A monoid $M$ is called an inverse monoid if, for every $a\in M$, there is a unique element $a^*\in M$ such that $aa^*a=a$ and $a^*aa^*=a^*.$ The idempotents of an inverse monoid commute. For a unit $a$, its inverse $a^*$ coincides with its group inverse $a^{-1}$; for an idempotent $e$, one has $e^*=e$.

Let $\IAM(D_n^\eps)$ denote the inverse Artin monoid of type $D_n$ defined in \cite[Definition~3.1]{Duan2025}. We use the vertex labelling of $D_n^\eps$ shown in Figure~\ref{the initial diagrams}: the vertices $1,\ldots,n$ are mutable, and $\eps$ is frozen. The mutable Dynkin diagram has edge set
\[
\mathcal{E}_D
=
\bigl\{\{1,3\},\{2,3\}\bigr\}
\cup
\bigl\{\{i,i+1\}:3\leq i<n\bigr\}.
\]

With this labelling, $\IAM(D_n^\eps)$ has generators $s_1^{\pm1},\ldots,s_n^{\pm1},s_\eps$ and defining relations
\[
\begin{aligned}
s_i s_i^{-1}=s_i^{-1}s_i&=1 && (1\leq i\leq n),\\
s_i s_j s_i&=s_j s_i s_j && (\{i,j\}\in\mathcal{E}_D),\\
s_i s_j&=s_j s_i && (i\ne j,\ \{i,j\}\notin\mathcal{E}_D),\\
s_\eps^2&=s_\eps,\\
s_\eps s_i&=s_i s_\eps && (1\leq i<n),\\
s_n^2 s_\eps&= s_\eps= s_\eps s_n^2,\\
s_\eps s_n s_\eps s_n&=s_n s_\eps s_n s_\eps= s_\eps s_n s_\eps,
\end{aligned}
\]
together with the fork relation
\[
s_1s_3\cdots s_n\, s_\eps\,s_n\cdots s_3s_1 = s_2s_3\cdots s_n\, s_\eps\,s_n\cdots s_3s_2.
\]
In particular, the generators $s_1,\ldots,s_n$ are units, with inverses $s_1^{-1},\ldots,s_n^{-1}$, and $s_\varepsilon$ is idempotent. Here we compose braids from right to left.

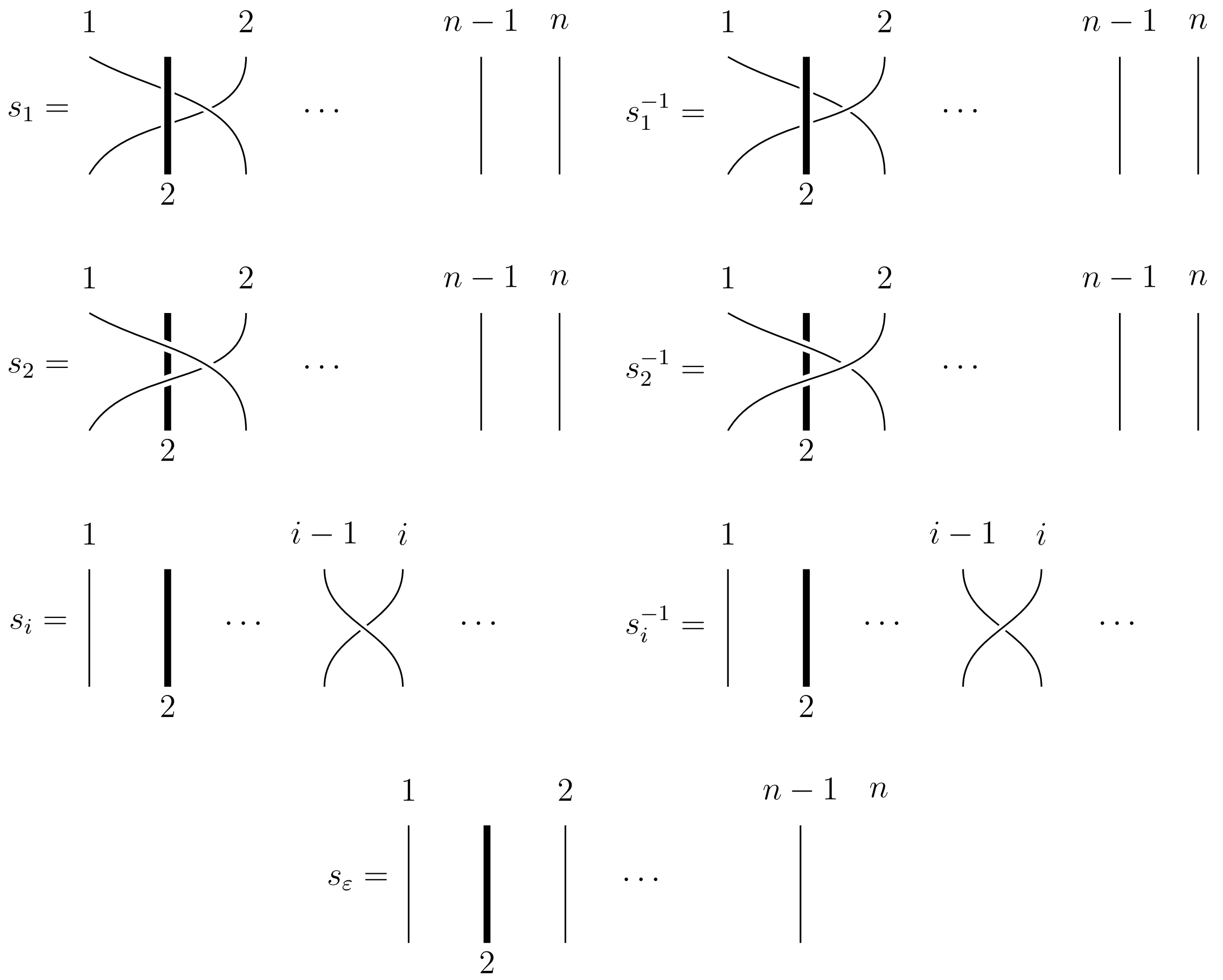
\begin{figure}[htbp]
\centering
\begin{tabular}{@{}c@{\hspace{0.03\textwidth}}c@{}}

\resizebox{0.46\textwidth}{!}{%
\begin{tikzpicture}[line width=0.6pt]
  \path[use as bounding box]
    (-1.35,-0.45) rectangle (6.3,2.15);

  \node at (-0.65,0.8) {$s_1=$};
  \node at (0,1.95) {$1$};
  \node at (2,1.95) {$2$};
  \node at (3,0.8) {$\cdots$};
  \node at (5,1.95) {$n-1$};
  \node at (6,1.95) {$n$};
  \node at (1,-0.25) {$2$};

  \draw (5,0)--(5,1.5);
  \draw (6,0)--(6,1.5);
  \draw (2,1.5) to[out=-90,in=60] (0,0);
  \draw[preaction={draw=white,line width=3pt}]
    (0,1.5) to[out=-30,in=90] (2,0);
  \draw[line width=2.5pt,
    preaction={draw=white,line width=5pt}]
    (1,0)--(1,1.5);
\end{tikzpicture}}
&
\resizebox{0.46\textwidth}{!}{%
\begin{tikzpicture}[line width=0.6pt]
  \path[use as bounding box]
    (-1.35,-0.45) rectangle (6.3,2.15);

  \node at (-0.8,0.8) {$s_1^{-1}=$};
  \node at (0,1.95) {$1$};
  \node at (2,1.95) {$2$};
  \node at (3,0.8) {$\cdots$};
  \node at (5,1.95) {$n-1$};
  \node at (6,1.95) {$n$};
  \node at (1,-0.25) {$2$};

  \draw (5,0)--(5,1.5);
  \draw (6,0)--(6,1.5);
  \draw (0,1.5) to[out=-30,in=90] (2,0);
  \draw[preaction={draw=white,line width=3pt}]
    (2,1.5) to[out=-90,in=60] (0,0);
  \draw[line width=2.5pt,
    preaction={draw=white,line width=5pt}]
    (1,0)--(1,1.5);
\end{tikzpicture}}
\\[1.2em]

\resizebox{0.46\textwidth}{!}{%
\begin{tikzpicture}[line width=0.6pt]
  \path[use as bounding box]
    (-1.35,-0.45) rectangle (6.3,2.15);

  \node at (-0.65,0.8) {$s_2=$};
  \node at (0,1.95) {$1$};
  \node at (2,1.95) {$2$};
  \node at (3,0.8) {$\cdots$};
  \node at (5,1.95) {$n-1$};
  \node at (6,1.95) {$n$};
  \node at (1,-0.25) {$2$};

  \draw (5,0)--(5,1.5);
  \draw (6,0)--(6,1.5);
  \draw[line width=2.5pt] (1,0)--(1,1.5);
  \draw[preaction={draw=white,line width=4pt}]
    (2,1.5) to[out=-90,in=60] (0,0);
  \draw[preaction={draw=white,line width=4pt}]
    (0,1.5) to[out=-30,in=90] (2,0);
\end{tikzpicture}}
&
\resizebox{0.46\textwidth}{!}{%
\begin{tikzpicture}[line width=0.6pt]
  \path[use as bounding box]
    (-1.35,-0.45) rectangle (6.3,2.15);

  \node at (-0.8,0.8) {$s_2^{-1}=$};
  \node at (0,1.95) {$1$};
  \node at (2,1.95) {$2$};
  \node at (3,0.8) {$\cdots$};
  \node at (5,1.95) {$n-1$};
  \node at (6,1.95) {$n$};
  \node at (1,-0.25) {$2$};

  \draw (5,0)--(5,1.5);
  \draw (6,0)--(6,1.5);
  \draw[line width=2.5pt] (1,0)--(1,1.5);
  \draw[preaction={draw=white,line width=4pt}]
    (0,1.5) to[out=-30,in=90] (2,0);
  \draw[preaction={draw=white,line width=4pt}]
    (2,1.5) to[out=-90,in=60] (0,0);
\end{tikzpicture}}
\\[1.2em]

\resizebox{0.46\textwidth}{!}{%
\begin{tikzpicture}[line width=0.6pt]
  \path[use as bounding box]
    (-1.35,-0.45) rectangle (6.3,2.15);

  \node at (-0.65,0.8) {$s_i=$};
  \node at (0,1.95) {$1$};
  \node at (1,-0.25) {$2$};
  \node at (2,0.8) {$\cdots$};
  \node at (3,1.95) {$i-1$};
  \node at (4,1.95) {$i$};
  \node at (5,0.8) {$\cdots$};

  \draw (0,0)--(0,1.5);
  \draw[line width=2.5pt] (1,0)--(1,1.5);
  \draw (4,1.5) to[out=-90,in=90] (3,0);
  \draw[preaction={draw=white,line width=3pt}]
    (3,1.5) to[out=-90,in=90] (4,0);
\end{tikzpicture}}
&
\resizebox{0.46\textwidth}{!}{%
\begin{tikzpicture}[line width=0.6pt]
  \path[use as bounding box]
    (-1.35,-0.45) rectangle (6.3,2.15);

  \node at (-0.8,0.8) {$s_i^{-1}=$};
  \node at (0,1.95) {$1$};
  \node at (1,-0.25) {$2$};
  \node at (2,0.8) {$\cdots$};
  \node at (3,1.95) {$i-1$};
  \node at (4,1.95) {$i$};
  \node at (5,0.8) {$\cdots$};

  \draw (0,0)--(0,1.5);
  \draw[line width=2.5pt] (1,0)--(1,1.5);
  \draw (3,1.5) to[out=-90,in=90] (4,0);
  \draw[preaction={draw=white,line width=3pt}]
    (4,1.5) to[out=-90,in=90] (3,0);
\end{tikzpicture}}
\\[1.2em]

\multicolumn{2}{c}{%
\resizebox{0.46\textwidth}{!}{%
\begin{tikzpicture}[line width=0.6pt]
  \path[use as bounding box]
    (-1.35,-0.45) rectangle (6.3,2.15);
  \node at (-0.65,0.8) {$s_\eps=$};
  \node at (0,1.95) {$1$};
  \node at (2,1.95) {$2$};
  \node at (3,0.8) {$\cdots$};
  \node at (5,1.95) {$n-1$};
  \node at (6,1.95) {$n$};
  \node at (1,-0.25) {$2$};

  \draw (0,0)--(0,1.5);
  \draw (2,0)--(2,1.5);
  \draw (5,0)--(5,1.5);
  \draw[line width=2.5pt] (1,0)--(1,1.5);
\end{tikzpicture}}}
\end{tabular}
\caption{Braid pictures of $s_1^{\pm1}$, $s_2^{\pm1}$, $s_i^{\pm1}$ for $3\leq i\leq n$, and the idempotent $s_\eps$. The thick vertical line represents the cone point of order $2$; the $n$th moving strand is absent in the picture for $s_\eps$.}\label{the braid picture of se D}
\end{figure}

Figure~\ref{the braid picture of se D} depicts braid representatives of $s_i^{\pm1}$, $1\le i\le n$, and $s_\varepsilon$. Following Allcock's conventions~\cite{Allcock2002}, the thick vertical line represents a cone point of order $2$, and the thin strands represent the moving strands, labelled $1,\ldots,n$. The idempotent $s_\varepsilon$ is represented by the partial braid obtained from the identity braid by deleting the $n$th moving strand.

\section{Inverse monoids associated with quivers in $\mathcal{M}(D_n^{\varepsilon})$}\label{monoids and diagrams}

We associate an inverse monoid $M_{\Gamma}$ with each quiver $\Gamma\in\mathcal{M}(D_n^{\varepsilon})$. Our main theorem asserts that its isomorphism type is invariant under mutation at mutable vertices.

\begin{definition}\label{def:B-inverse-monoid}
Let $\Gamma\in\mathcal{M}(D_n^{\varepsilon})$ have mutable vertex set $I=\{1,\ldots,n\}$, with $n\geq 4$, and frozen vertex $\varepsilon$. Define $M_{\Gamma}$ to be the inverse monoid with generators
\[
s_i^{\pm1}\quad(i\in I),\qquad s_{\varepsilon},
\]
subject to the following relations.

\begin{itemize}
\item[(R1)] For distinct $i,j\in I$,
\[
\begin{aligned}
s_is_j&=s_js_i
&&\text{if $i$ and $j$ are not adjacent},\\
s_is_js_i&=s_js_is_j
&&\text{if $i$ and $j$ are adjacent}.
\end{aligned}
\]

\item[(R2)] $s_is_i^{-1}=s_i^{-1}s_i=1$ for $i\in I$, and $s_\eps^2=s_\eps$.

\item[(R3)] $s_\eps s_i=s_is_\eps$ when $i$ is not adjacent to $\eps$.

\item[(R4)] When $i$ is adjacent to $\eps$,
\[
s_is_\eps s_is_\eps=s_\eps s_is_\eps s_i=s_\eps s_is_\eps,
\qquad s_i^2s_\eps=s_\eps=s_\eps s_i^2.
\]

\item[(R5)] For every entirely mutable oriented chordless cycle $C$
\[
C:\quad i_0\longrightarrow i_1\longrightarrow\cdots \longrightarrow i_{d-1}\longrightarrow i_0,
\]
we impose
\[
[s_{i_a},p_a(C;s)]=1 \qquad(0\leq a<d).
\]

\item[(R6)] For every oriented triangle $\varepsilon\longrightarrow i\longrightarrow j\longrightarrow\varepsilon$, we impose
\[
s_{\varepsilon}s_js_is_j^{-1} = s_js_is_j^{-1}s_{\varepsilon}.
\]

\item[(R7)]  All paths used here lie in underlying undirected graphs. For a shortest path $\pi=(v_1,\ldots,v_r,\eps)$, let
\[
P_\pi=s_{v_1}\cdots s_{v_r}s_\eps,
\qquad \overline P_\pi=s_\eps s_{v_r}\cdots s_{v_1}.
\]
For the trivial path $(\eps)$ set $P_\pi=\overline P_\pi=s_\eps$. In Types I--III, take shortest paths $\pi_a,\pi_b$ from $a,b$ to $\eps$ and impose
\[P_{\pi_a}\overline P_{\pi_a}=P_{\pi_b}\overline P_{\pi_b}.\]
In Type IV, take a shortest path $\pi_{c'}$ from $c'$ to $\eps$ within $Q'$, set $E_{c'}=P_{\pi_{c'}}\overline P_{\pi_{c'}}$, and impose
\[s_{i_1}E_{c'}s_{i_1}=s_{i_2}\cdots s_{i_d}E_{c'}s_{i_d}\cdots s_{i_2}.\]
If $c'=\eps$, then $E_{c'}=s_\eps$.
\end{itemize}
\end{definition}

Note that the shortest paths used in (R7) are unambiguously determined in Types~I--IV, since shortest paths within type~$A$ components are unique. For $\Gamma\in\mathcal{M}(A_n^{\varepsilon})$, the monoid $M_{\Gamma}$ defined by the same generators subject to relations (R1)--(R6)
coincides with the monoid defined in \cite[Definition 4.1]{Duan2025}.

\begin{theorem}\label{thm:B-mutation-invariance}
Let $\Gamma\in \mathcal{M}(D_n^{\varepsilon})$ and let $k$ be mutable. The assignment
\[
\Phi_k:M_{\mu_k(\Gamma)}\longrightarrow M_\Gamma,
\qquad t_v \longmapsto 
\begin{cases}s_ks_vs_k^{-1},&v\to k\text{ in }\Gamma,\\s_v,&\text{otherwise},\end{cases}
\]
for $v\in I\sqcup\{\eps\}$ is an isomorphism. Consequently,
\[M_\Gamma\cong\IAM(D_n^{\varepsilon}).\]
\end{theorem}
For an initial Dynkin orientation, (R5) and (R6) are absent and (R7) is the standard fork relation. The last assertion in Theorem \ref{thm:B-mutation-invariance} therefore follows from mutation invariance and the presentation in Section \ref{sec:nverse Artin monoid of type Dn}. Thus each quiver in this mutation class gives a presentation of the inverse Artin monoid $\IAM(D_n^{\varepsilon})$.

\section{Proof of Theorem~\ref{thm:B-mutation-invariance}} \label{proof of main theorem}

\subsection{Inverse half twists}\label{sec:inverse half twists}

Let $X$ be an oriented disk with a finite set $P$ of marked points in its interior. Let $D\subset\operatorname{int}(X)$ be an embedded closed disk such that
\[
D\cap P=\{A,B,C\}\subset\partial D,
\]
where $A,B,C$ are distinct and occur in clockwise order on $\partial D$. Denote the three successive boundary arcs of $D$ by $AB,BC,CA$.

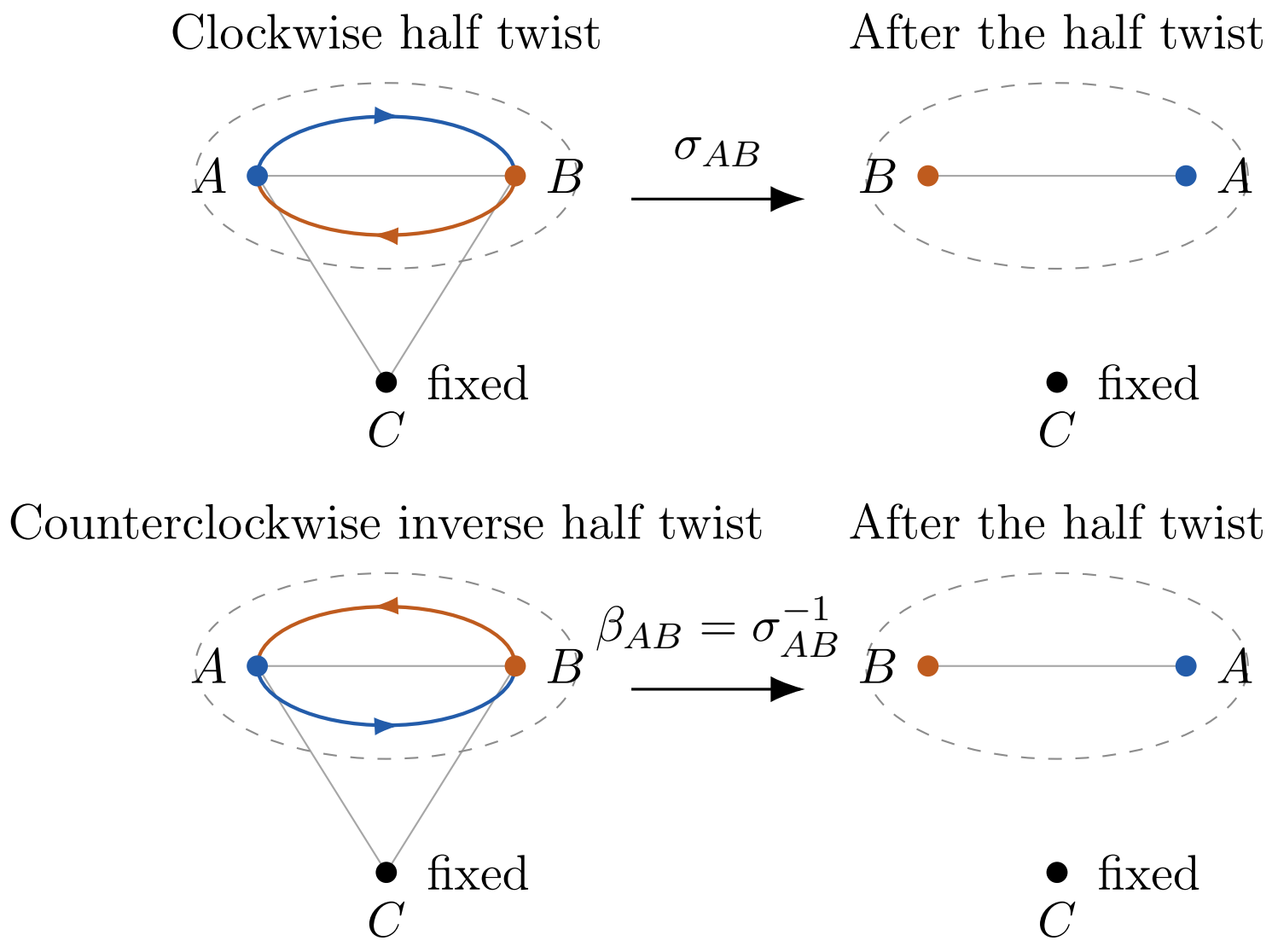
\begin{figure}[H]
\begin{tikzpicture}[font=\small]
\TwistRow{0cm}{1}{\sigma_{AB}}{Clockwise half twist}
\TwistRow{-3.8cm}{-1}{\beta_{AB}=\sigma_{AB}^{-1}}{Counterclockwise inverse half twist}
\node[align=center] at (2.6,-6.25){$A\longmapsto B,\qquad B\longmapsto A,\qquad C\longmapsto C.$};
\end{tikzpicture}
\caption{A clockwise half twist and its inverse, taken in the ambient disk $X$.}\label{half twist}
\end{figure}

Following the convention of \cite[Definition 3.3]{GM2017}, let $\sigma_{AB}$ be the clockwise half twist along $AB$, taken in the ambient disk $X$, and set
\[
\beta_{AB}=\sigma_{AB}^{-1}.
\]
Thus $\beta_{AB}$ exchanges $A$ and $B$ counterclockwise in a sufficiently small neighbourhood of $AB$, while fixing every point of $P\setminus\{A,B\}$, see Figure~\ref{half twist}. Define $\sigma_{BC},\sigma_{CA}$ and their inverses $\beta_{BC},\beta_{CA}$ similarly.

By \cite[Lemma 3.5]{GM2017} and cyclic relabelling,
\[
\sigma_{AB}\sigma_{BC}
=\sigma_{BC}\sigma_{CA}
=\sigma_{CA}\sigma_{AB}.
\]
Taking inverses gives the corresponding relations for inverse half twists:
\[
\beta_{BC}\beta_{AB}
=\beta_{CA}\beta_{BC}
=\beta_{AB}\beta_{CA}.
\]
Consequently,
\[
\beta_{BC}\beta_{AB}\beta_{BC}^{-1}=\beta_{CA}.
\]
Hence conjugation by \(\beta_{BC}\) transforms the inverse half twist along \(AB\) into the inverse half twist along \(CA\), carrying its endpoints from \(\{A,B\}\) to \(\{A,C\}\). These relations also imply the usual braid relation
\[
\beta_{AB}\beta_{BC}\beta_{AB}
=\beta_{BC}\beta_{AB}\beta_{BC}.
\]

\subsection{A geometric realization of the inverse Artin monoid of type
$A_m^\varepsilon$}\label{sec:A-geometric-preliminaries}

Let $Q\in\mathcal M(A_m^\varepsilon)$. Choose a triangulation $\mathcal T$ of an $(m+3)$-gon whose internal diagonals are labelled by the mutable
vertices and whose distinguished boundary side is labelled $\varepsilon$, so that its associated labelled ice quiver is $Q$. We write $T=T_Q^*$ for the embedded dual tree: its vertices lie in the triangles, its edges cross the internal diagonals, and its root $o$ lies in the triangle incident with $\varepsilon$. Thus $P=V(T)$ has $m+1$ vertices and $T$ has $m$ edges. The root has degree at most two, and a mutable vertex is adjacent to $\varepsilon$ precisely when its corresponding edge is incident with $o$.

Let $B(P)$ and $\IB(P)$ denote the braid group and the inverse braid monoid of the ambient disk with marked point set $P$, respectively. For $x\in P$, let $\delta_x$ be the idempotent obtained from the identity braid by deleting the strand at $x$. If $\beta$ is an inverse half twist with endpoints $u,v$ and $\tau=(u\ v)$, the deletion identities are
\begin{align}
 \beta\delta_x\beta^{-1}
 &=\delta_{\tau(x)}
 &&(x\in P),\label{delconj}\\
 \beta^2\delta_x
 &=\delta_x=\delta_x\beta^2
 &&(x=u,v),\label{delsquare}\\
 \beta\delta_u\delta_v
 &=\delta_u\delta_v=\delta_u\delta_v\beta.
 &&\label{deldouble}
\end{align}
See \cite[Theorem~2.1]{Vershinin2009}. For either ordering $(x,y)$ of $\{u,v\}$, these imply
\begin{align}
 \beta\delta_x\beta
 &=\delta_y,\label{positiveconj}\\
 \beta\delta_x\beta\delta_x
 &=\delta_x\beta\delta_x\beta
 =\delta_x\beta\delta_x
 =\delta_x\delta_y.\label{mixedlocal}
\end{align}
Indeed, \eqref{delconj} and \eqref{delsquare} give $\beta\delta_x\beta=\delta_y\beta^2=\delta_y$, while \eqref{delconj} and \eqref{deldouble} give $\delta_x\beta\delta_x=\delta_x\delta_y\beta=\delta_x\delta_y$.

\begin{proposition}\label{atype}
Let $U_Q$ be the inverse monoid defined by (R1)--(R6) for $Q$. For each mutable vertex $i$, let $\beta_i$ be the inverse half twist along edge $i$ of $T$. Then the assignments
\[
\Theta_Q:U_Q\longrightarrow\IB(P),\qquad s_i\longmapsto\beta_i,\qquad s_\varepsilon\longmapsto\delta_o
\]
define an isomorphism of inverse monoids.
\end{proposition}

\begin{proof}
The case $m=0$ is immediate. Assume $m\geq1$. By \cite[Theorem~3.6]{GM2017} and \cite[Definitions~4.4--4.5 and Theorem~6.1]{HHLP2017}, using inverse half twists to match our conventions, the mutable relations (R1) and (R5) present $B(P)$ on the generators $\beta_i$. The relations (R2)--(R4) follow from the idempotence of $\delta_o$ and \eqref{delconj}, \eqref{delsquare}, and \eqref{mixedlocal}. For a frozen triangle
$\varepsilon\to i\to j\to\varepsilon$, write the endpoints of $i,j$ as $\{u,o\}$ and $\{o,v\}$, respectively. The braid $\beta_j\beta_i\beta_j^{-1}$ induces the transposition $(u\ v)$ and hence fixes $o$. Repeated application of \eqref{delconj} therefore gives
\[
\delta_o\beta_j\beta_i\beta_j^{-1}=\beta_j\beta_i\beta_j^{-1}\delta_o,
\]
which is (R6). Thus $\Theta_Q$ is a homomorphism.

To construct the inverse, the mutable presentation above defines a group homomorphism $\Lambda:B(P)\to U_Q^\times$ by $\Lambda(\beta_i)=s_i$, since the units $s_i$ satisfy (R1) and (R5). We verify that it extends by $\Lambda(\delta_o)=s_\varepsilon$. Choose an edge $a$ of $T$ incident with $o$ and a closed disk $D$ with $D\cap P=P\setminus\{o\}$, meeting $a$ only in a terminal subarc, as follows.
If $o$ has degree one, take $a$ to be its unique incident edge and
$D$ a small regular neighbourhood of the tree obtained by deleting
$o$ and $a$.
If $o$ has degree two, label its incident edges $i,j$ so that
$\varepsilon\to i\to j\to\varepsilon$, and write
$i=\{u,o\}$, $j=\{o,v\}$.
Choose an arc $c$, with interior disjoint from $T$, such that
$i,j,c$ bound a thin disk with $u,o,v$ in clockwise order.
Such an arc is obtained by taking a sufficiently thin disk on the
appropriate side of $u-o-v$, narrowing it near $u$ and $v$ to avoid
the other incident edges; no third edge meets $o$.
By \cite[Lemma~3.5]{GM2017},
\[
 \beta_c=\beta_j\beta_i\beta_j^{-1}.
\]
Take $a=j$ and $D$ a small regular neighbourhood of the tree obtained
by deleting $o,i,j$ and adding $c$.

By \cite{Sergiescu1993,Sergiescu1994}, $B(D,P\setminus\{o\})$ is generated by the unchanged edge half twists and, in the second case, $\beta_c$. Their images under $\Lambda$ commute with $s_\varepsilon$ by (R3) and, for $\Lambda(\beta_c)=s_js_is_j^{-1}$, by (R6). Hence $s_\varepsilon$ commutes with the image of this entire subgroup.

Put $D$ in standard position and choose an Artin chain starting at $o$, whose first edge is $a$ and whose remaining edges lie in $D$. The images of the latter half twists commute with $s_\varepsilon$, while the first has image $s_a$. Thus the idempotence relation is (R2), the required commutation relations hold by the preceding paragraph, and the square and mixed relations with the first generator are (R4). Together with the Artin relations, these are the relations of the standard inverse braid monoid presentation \cite[relations~(1.7)]{Vershinin2009}. This presentation also holds for inverse half twists: (R4) gives
\[
s_a^{-1}s_\varepsilon=s_as_\varepsilon, \qquad  s_\varepsilon s_a^{-1}=s_\varepsilon s_a,
\]
so the deletion relations are unchanged when the Artin generators are inverted. Consequently, $\Lambda$ extends to a homomorphism $\IB(P)\to U_Q$ with $\Lambda(\delta_o)=s_\varepsilon$.

Hence $\Theta_Q$ and $\Lambda$ are mutually inverse.
\end{proof}

\begin{lemma}\label{transport}
Suppose elements $s_1,\ldots,s_m,s_\eps$ in an inverse monoid $N$ satisfy (R1)--(R6) for $Q$. There are idempotents $F_x\in N$, indexed by $V(T_Q^*)$, with $F_o=s_\eps$, such that:
\begin{enumerate}[label=(\roman*)]
\item If edge $i$ has endpoints $u_i,v_i$ and $\tau_i=(u_i\ v_i)$, then
\[s_iF_xs_i^{-1}=F_{\tau_i(x)}.\]
\item For either ordering $(x,y)$ of these endpoints,
\[
s_i^2F_x=F_x=F_xs_i^2,
\qquad
s_iF_xs_iF_x=F_xs_iF_xs_i=F_xs_iF_x=F_xF_y.
\]
\item If $\pi=(j_1,\ldots,j_r,\eps)$ is a shortest path in $Q$, then the dual edges $j_r,\ldots,j_1$ form a path from $o$ to a vertex $x_\pi$, and
\[
s_{j_1}\cdots s_{j_r}s_\eps s_{j_r}\cdots s_{j_1}
=w_\pi s_\eps w_\pi^{-1}=F_{x_\pi},
\qquad w_\pi=s_{j_1}\cdots s_{j_r}.
\]
The case $r=0$ gives $F_o=s_\eps$.
\end{enumerate}
\end{lemma}
\begin{proof}
The homomorphism $U_Q\to N$ and Proposition~\ref{atype} give a homomorphism $\rho:\IB(P)\to N$. Set $F_x=\rho(\delta_x)$. Parts (i) and (ii) follow from \eqref{delconj}, \eqref{delsquare}, and \eqref{mixedlocal}.

For (iii), assume first that $r\ge1$. Adjoin to $T_Q^*$ a new leaf $o_\varepsilon$ and a formal edge labelled $\varepsilon$ joining $o_\varepsilon$ to the root $o$. Denote the resulting tree by $\widehat T_Q^*$. The underlying undirected graph of $Q$ is the line graph of $\widehat T_Q^*$.

The reversed sequence $(\varepsilon,j_r,j_{r-1},\ldots,j_1)$ is therefore a shortest path in this line graph. Its edge labels are distinct. If three consecutive edges met at the same vertex of $\widehat T_Q^*$, the first and third would be adjacent in the line graph, giving a shortcut. Thus the consecutive intersection vertices are distinct.
Since $\widehat T_Q^*$ is a tree, the displayed edge sequence follows a simple tree path starting at $o_\varepsilon$. After its first edge $\varepsilon$, this path starts at $o$ and traverses $j_r,\ldots,j_1$. Let $x_\pi$ be its terminal vertex.

If an edge labelled $i$ joins successive vertices $x,y$ on this tree path, then parts (i) and (ii) give
\[
s_iF_xs_i^{-1}=F_y, \qquad  F_xs_i^2=F_x.
\]
The second identity implies $F_xs_i=F_xs_i^{-1}$, and hence $s_iF_xs_i=s_iF_xs_i^{-1}=F_y.$  Iterating these identities along the path, starting with $F_o=s_\varepsilon$, yields
\[
\begin{aligned}
F_{x_\pi} &=w_\pi s_\varepsilon w_\pi^{-1} =s_{j_1}\cdots s_{j_r}\,
s_\varepsilon\, s_{j_r}\cdots s_{j_1}, \qquad  w_\pi=s_{j_1}\cdots s_{j_r}.
\end{aligned}
\]
For $r=0$, take $x_\pi=o$ and $w_\pi=1$; the assertion reduces to $F_o=s_\varepsilon$.
\end{proof}

\subsection{Idempotent separation and the proof of Theorem~\ref{thm:B-mutation-invariance}}

Let $\Bool_\Gamma$ be the inverse monoid defined by relations (R1)--(R7) for $\Gamma$, together with the additional relations $s_i^2=1$ for all mutable vertices $i$. By definition, there is a canonical epimorphism
\[
q = q_\Gamma:M_\Gamma\twoheadrightarrow\Bool_\Gamma.
\]
By \cite[Definition 4.5 and Theorems 4.7 and 4.9]{DLL2020}, $\Bool_\Gamma$ is isomorphic to the Boolean reflection monoid of type $D_n$. Indeed, let $r_i=q_\Gamma(s_i)$ and $e=q_\Gamma(s_\varepsilon)$. Since $r_i^{-1}=r_i$, relations (R1)--(R5) specialize to \cite[Definition 4.5, (R1), (R2), and (R3)(i)]{DLL2020}. The braid relation identifies (R6) with the frozen triangle relation (R3)(ii) there. For $\pi=(v_1,\ldots,v_r,\varepsilon)$, let $w_\pi=r_{v_1}\cdots r_{v_r}$. Then
\[
q_\Gamma(P_\pi\overline P_\pi)=w_\pi e w_\pi^{-1},
\]
so (R7) is precisely \cite[Definition 4.5, (R4)(ii)]{DLL2020}. For the frozen four-cycle $\varepsilon\to a\to d\to b\to\varepsilon$, (R7) gives $r_aer_a=r_ber_b$, hence $er_ar_b=r_ar_be$. Since $er_d=r_de$ and
\[
r_ar_dr_br_dr_a=(r_ar_b)r_d(r_ar_b),
\]
the required frozen four-cycle relation follows. Thus the presentations coincide, and \cite[Theorems 4.7 and 4.9]{DLL2020} give $B_\Gamma\cong \mathcal B(D_n)$.

\begin{lemma}\label{forest}
In Types~I--III, let $Q_L=\Gamma\setminus\{b\}$, $Q_R=\Gamma\setminus\{a\}.$  In Type~IV, let
\[
Q_L=\Gamma\setminus\{i_d\}, \qquad Q_R=\Gamma\setminus\{i_1\}.
\]
Then $Q_L,Q_R\in\mathcal M(A_{n-1}^{\varepsilon})$. Their rooted dual trees can be identified on a common vertex set so that deleting edge $a$ from $T_{Q_L}^*$ and edge $b$ from $T_{Q_R}^*$ in Types~I--III, or edge $i_1$ from $T_{Q_L}^*$ and edge $i_d$ from $T_{Q_R}^*$ in Type~IV, leaves the same edge-labelled forest
\[
F=F_0\sqcup F_1,
\]
with a common root $o\in F_0$. All $n$ vertices are retained. Moreover, the root-to-vertex paths determined by the two sides of {\rm(R7)} end at the same vertex of $F_1$.
\end{lemma}

\begin{proof}
We first specify the convention for free endpoints.
For a rooted dual tree, adjoin a formal edge labelled
$\varepsilon$ joining the root to a new leaf.
A free endpoint of a mutable edge means a leaf endpoint
in this augmented tree. For a type~$A$ component without
a frozen vertex, it means a leaf endpoint in its ordinary
dual tree. The edge corresponding to a mutable connecting
vertex has such an endpoint, by the line-graph description
of type~$A$ quivers. The formal edge and its added vertex
are used only to specify the gluing positions and are
omitted from all forests below.

The descriptions of Types~I--IV show that $Q_L,Q_R\in\mathcal M(A_{n-1}^{\varepsilon})$. We construct the required identifications case by case.

\textbf{Type I.} Let $F_0=T_{Q'}^*$, with its root $o$, and choose
the free endpoint $u$ of the edge labelled $c$.
Let $F_1=\{v\}$ consist of one new vertex.
The tree $T_{Q_L}^*$ is obtained by joining $u$ to $v$
with edge $a$, whereas $T_{Q_R}^*$ is obtained by joining
the same vertices with edge $b$.
Deleting these edges leaves $F_0\sqcup F_1$.
The two root paths in {\rm(R7)} follow the same path
from $o$ to $u$ and then cross $a$ or $b$, respectively.
Both end at $v$.

\textbf{Type II.} Take the dual trees $T_{Q'}^*$ and $T_{Q''}^*$. If $c\ne\varepsilon$, choose the free endpoint of
edge $c$ in $T_{Q'}^*$; if $c=\varepsilon$, use the
root $o$ of $T_{Q'}^*$ instead.
Identify this vertex with a free endpoint of edge $d$
in $T_{Q''}^*$, and denote the identified vertex by $u$.
Let $F_0$ be the resulting tree, rooted at the image
of $o$, and let $F_1=\{v\}$ consist of one new vertex. Adjoining an edge from $u$ to $v$, labelled $a$ or $b$,
gives $T_{Q_L}^*$ or $T_{Q_R}^*$, respectively.
These constructions realize the triangles with labels
$\{a,c,d\}$ and $\{b,c,d\}$.
When $c=\varepsilon$, the formal frozen edge occupies
the position of $c$ in this description.
Deleting the added edges leaves $F_0\sqcup F_1$.
Both root paths in {\rm(R7)} end at $v$.

\textbf{Type III.} Let $F_0=T_{Q'}^*$, rooted at $o$, and let $F_1=T_{Q''}^*$. If $c\ne\varepsilon$, let $u\in F_0$ be the free endpoint of edge $c$; otherwise, let $u=o$. Choose a free endpoint $v\in F_1$ of edge $d$. The trees $T_{Q_L}^*$ and $T_{Q_R}^*$ are obtained from $F_0\sqcup F_1$ by joining $u$ to $v$ with edge $a$ and edge $b$, respectively. Deleting these edges therefore gives the same forest. Both root paths in {\rm(R7)} end at $v$.

%Use the numbering
%\[
%    i_1\to i_2\to\cdots\to i_d\to i_1,
%    \qquad
%    i_d\to i_1\to c'\to i_d.
%\]
%
\textbf{Type IV.}  Let $F_0=T_{Q'}^*$, rooted at $o$. If $c'\ne\varepsilon$, let $u$ be the free endpoint of edge $c'$; otherwise, set $u=o$. Start with the path on vertices $v_1,\ldots,v_{d-1}$, where edge $i_j$ joins $v_{j-1}$ to $v_j$ for $2\le j\le d-1$. For each spike other than the one with tip $c'$, say on $i_j\to i_{j+1}$ with $1\le j\le d-1$, attach the dual tree of its component by identifying a free endpoint of its tip edge with $v_j$. Denote the resulting tree by $F_1$. The tree $T_{Q_L}^*$ is obtained by adding edge $i_1$ between $u$ and $v_1$, whereas $T_{Q_R}^*$ is obtained by adding edge $i_d$ between $u$ and $v_{d-1}$. Deleting these added edges leaves the same forest $F_0\sqcup F_1$, with all other edge labels unchanged. After following the common path from $o$ to $u$, the two dual paths determined by {\rm(R7)} traverse the edge sequences
\[
(i_1) \qquad\text{and}\qquad  (i_d,i_{d-1},\ldots,i_2),
\]
respectively. Both end at $v_1$.

In every case, the root lies in the same component $F_0$. Each of the two dual trees has $n$ vertices. Only the indicated mutable edges are deleted, so all $n$ vertices are retained.
\end{proof}

We now identify $B_\Gamma$ with the Boolean reflection monoid $\mathcal B(D_n)$ in its standard realization on $\mathbb R^n$, using the presentation isomorphism established above. We choose this identification by transporting the standard realization for an initial Dynkin orientation along the Boolean mutation isomorphisms in \cite[Theorem 4.7 and Proposition 7.1]{DLL2020}. In this realization,
\[
\mathcal B(D_n)= \left\{ w|_U : w\in W(D_n),\  U\text{ is a coordinate subspace of }\mathbb R^n \right\},
\]
where $W(D_n)$ acts by signed permutation matrices with an even number of negative entries. We continue to write $q$ for the quotient map composed with this identification.

For $1\le x\le n$, let
\[
V_x=\{a=(a_1,\ldots,a_n)\in\mathbb R^n:a_x=0\}, \qquad \delta_x=\operatorname{id}_{V_x}.
\]
In what follows, $\delta_x$ denotes this partial linear identity in $\mathcal B(D_n)$. Let $H$ be the subgroup of units of $M_\Gamma$ generated by $s_i$, $i\in I$, and let $\mathcal C=\{hs_\varepsilon h^{-1}:h\in H\}.$  Then $q(H)=W(D_n)$, each $q(s_i)$ is a reflection, and $q(s_\varepsilon)\in\{\delta_1,\ldots,\delta_n\}.$

\begin{proposition}\label{orbit}
The conjugacy orbit $C$ consists of exactly $n$ idempotents $E_1,\ldots,E_n$, which can be indexed so that $q(E_x)=\delta_x$. Every $H$-conjugate of a mutable generator acts on $C$ by a transposition. For such an element $u$ and any $f\in C$, one has
\[
ufu^{-1}\ne f \quad\Longleftrightarrow\quad  q(u)q(f)q(u)q(f)=q(f)q(u)q(f).
\]
Whenever these equivalent conditions hold,
\begin{equation}\label{localD}
u^2f=f=fu^2, \qquad  ufuf=fufu=fuf.
\end{equation}
\end{proposition}

\begin{proof}
Apply Lemma~\ref{transport} to $Q_L$ and $Q_R$ inside $M_\Gamma$, and index the resulting families $F_z^L,F_z^R$ by the common forest in Lemma~\ref{forest}. They agree at the root, so transport along the common edges gives equality on $F_0$. By Lemma~\ref{forest}, the two paths in (R7) end at the same vertex of $F_1$. Thus (R7) and Lemma~\ref{transport}(iii) give equality of the two families at that vertex. Transport along the common edges then gives equality throughout $F_1$. Set $E_z=F_z^L=F_z^R$.

Every mutable generator occurs in at least one of $Q_L,Q_R$. Hence Lemma~\ref{transport}(i) shows that the set of the $E_z$ is invariant under $H$. It contains $s_\eps$, and every $E_z$ is obtained from $s_\eps$ by transport along a path in either dual tree. Consequently this set is exactly $C$, and $|C|\le n$. Since $q(H)=W(D_n)$, the image $q(C)$ is the full orbit of a coordinate-hyperplane idempotent. Therefore
\[
q(C)=\{\delta_1,\ldots,\delta_n\}, \qquad n=|q(C)|\le |C|\le n.
\]
Thus the $E_z$ are distinct and can be indexed as stated.

By Lemma~\ref{transport}(i), each mutable generator acts on $C$ by its endpoint transposition. Since the $E_z$ are distinct, Lemma~\ref{transport}(ii) gives the local identities in \eqref{localD} for $u=s_i$ and each idempotent moved by $s_i$. Now let $u=vs_iv^{-1}$ with $v\in H$. Its action on $C$ is conjugate to the action of $s_i$, so it is again a transposition. If $ufu^{-1}\ne f$, then $v^{-1}fv\in C$ is moved by $s_i$. Applying the preceding local identities to $s_i$ and $v^{-1}fv$, and conjugating by $v$, proves \eqref{localD}. Taking Boolean images gives the forward implication in the displayed equivalence.

Conversely, suppose that the Boolean relation holds but $ufu^{-1}=f$. Then $q(u)$ and $q(f)$ commute. Since $q(u)^2=1$ and $q(f)^2=q(f)$, the Boolean relation reduces to
\[
q(f)=q(u)q(f).
\]
Thus $q(u)$ fixes the coordinate hyperplane supporting $q(f)$ pointwise. This is impossible: $q(u)$ is a type $D_n$ reflection, which exchanges two coordinate axes, possibly changing both signs, and every coordinate hyperplane contains at least one of these axes. Hence $ufu^{-1}\ne f$, proving the reverse implication.
\end{proof}

\begin{proposition}\label{idemsep}
Every idempotent of $M_\Gamma$ has a unique expression
\[
E_S=\prod_{x\in S}E_x\qquad(S\subseteq\{1,\ldots,n\}).\]
In particular, $q$ is idempotent-separating and $M_\Gamma$ has exactly $2^n$ idempotents.
\end{proposition}

\begin{proof}
All elements of $C$ are idempotent and hence commute. We claim that $M_\Gamma$ is factorizable. More precisely, every element of $M_\Gamma$ has the form $Ph$, where $P$ is a product of conjugates of $s_\eps$ and $h\in H$. Indeed, right multiplication by a unit preserves this form, whereas
\[
(Ph)s_\eps=P(hs_\eps h^{-1})h.
\]
Since $M_\Gamma$ is generated by $H\cup\{s_\eps\}$, the assertion follows by induction on word length, starting with $P=h=1$. As $P$ is a product of commuting idempotents, we have $P^2=P$ and $P^*=P$.

Let $f\in M_\Gamma$ be idempotent and write $f=Ph$.
Then
\[
f^*=(Ph)^*=h^{-1}P,
\]
and, since $f^*=f$ and $f^2=f$,
\[
f=ff^*=(Ph)(h^{-1}P)=P^2=P.
\]
By Proposition~\ref{orbit}, every factor of $P$ is one of $E_1,\ldots,E_n$. Their commutativity and idempotence allow us to remove repeated factors, giving
\[
f=E_S=\prod_{x\in S}E_x
\]
for some $S\subseteq\{1,\ldots,n\}$. Here $E_\varnothing=1$. Conversely, every such product is idempotent.

Since $q(E_x)=\delta_x$, we have
\[
q(E_S) =\prod_{x\in S}\delta_x =\operatorname{id}_{\bigcap_{x\in S}V_x} =\operatorname{id}_{\{a\in\mathbb R^n:a_x=0\text{ for every }x\in S\}}.
\]
Here the intersection over the empty set is understood to be $\mathbb R^n$.

If $S\ne T$, choose $j\in S\mathbin{\triangle}T$. The $j$th standard basis vector belongs to exactly one of the two domains, so $q(E_S)\ne q(E_T)$. Thus the subset $S$ representing an idempotent is unique. The same argument shows that two idempotents with the same image under $q$ are equal. Hence $q$ is idempotent-separating. Finally, the bijection $S\mapsto E_S$ shows that $M_\Gamma$ has exactly $2^n$ idempotents.
\end{proof}

\begin{corollary}\label{lift}
If $u$ is a unit and $f$ an idempotent of $M_\Gamma$, then
\[
q(u)q(f)=q(f)q(u)\quad\Longrightarrow\quad uf=fu.
\]
Two idempotents with the same Boolean image are equal.
\end{corollary}
\begin{proof}
Since $ufu^{-1}$ and $f$ are idempotents with the same image under $q$, Proposition~\ref{idemsep} gives $ufu^{-1}=f$. Multiplying on the right by $u$ yields $uf=fu$.
\end{proof}

\begin{proof}[Proof of Theorem~\ref{thm:B-mutation-invariance}]
Let $\Gamma'=\mu_k(\Gamma)$, and let $s_v$ and $t_v$ be the generators of $M_\Gamma$ and $M_{\Gamma'}$, respectively. Here $q:M_\Gamma\to \Bool_\Gamma (\cong \Bool(D_n))$ is the epimorphism to the associated Boolean reflection monoid of type $D_n$.
Set
\[
r_v=
\begin{cases}
s_ks_vs_k^{-1},&v\to k\text{ in }\Gamma,\\
s_v,&\text{otherwise}.
\end{cases}
\]
We verify that the elements $r_v$ satisfy (R1)--(R7) for $\Gamma'$.

Deleting $\eps$ commutes with mutation at $k$. Theorem~\ref{thm:mutable} therefore gives (R1) and (R5) for the mutable elements $r_i$. Each $r_i$ is a unit, and $r_\eps$ is a conjugate of $s_\eps$, so (R2) also holds.

By the Boolean mutation theorem \cite[Theorem 4.7 and Proposition 7.1]{DLL2020}, the elements $q(r_v)$ are the Boolean generators associated
to $\Gamma'$. If $i$ is not adjacent to $\eps$ in $\Gamma'$, then $q(r_i)$ and $q(r_\eps)$ commute. Corollary~\ref{lift}, applied to $r_i$ and $r_\eps$, gives (R3). For an oriented triangle $\eps\to i\to j\to\eps$, the Boolean frozen cycle relation, together with the braid relation, gives
\[
q(r_\eps)q(r_jr_ir_j^{-1})
=q(r_jr_ir_j^{-1})q(r_\eps).
\]
Applying Corollary~\ref{lift} to $r_jr_ir_j^{-1}$ and $r_\eps$ gives (R6).

For (R4), suppose that $i$ is adjacent to $\eps$ in $\Gamma'$. The element $r_i$ is an $H$-conjugate of a mutable generator, and $r_\eps\in C$.
Their Boolean images under $q$ satisfy
\[
q(r_i)q(r_\eps)q(r_i)q(r_\eps)
=q(r_\eps)q(r_i)q(r_\eps).
\]
Proposition~\ref{orbit}, applied to $r_i$ and $r_\eps$, therefore gives (R4).

For (R7), we choose the two type $A$ subquivers of $\Gamma'$ from Lemma~\ref{forest}. The corresponding elements $r_v$ satisfy (R1)--(R6), as
proved above. In Types I--III, the two paths in (R7) start at $a$ or $b$ and then follow the shortest route through $c$ to $\eps$. In Type IV, they follow $i_1$ or $i_2,\ldots,i_d$, respectively, and then the shortest route from $c'$ to $\eps$ inside the attached component. The tails are trivial when $c=\eps$ or $c'=\eps$. These are shortest paths in the respective type $A$ subquivers. Using $r_\eps^2=r_\eps$, Lemma~\ref{transport}(iii) shows that both sides of (R7), evaluated at the elements $r_v$, are idempotents, and their images under $q$ are the same, so Proposition~\ref{idemsep} shows that they are equal.

Hence the assignment $t_v\mapsto r_v$ defines a homomorphism $\Phi_k:M_{\Gamma'}\longrightarrow M_\Gamma.$ Applying the same argument to the reverse mutation gives a homomorphism
\[
\widetilde\Phi_k:M_\Gamma\longrightarrow M_{\Gamma'}, \qquad s_v\longmapsto
\begin{cases}
t_kt_vt_k^{-1},&v\to k\text{ in }\Gamma',\\
t_v,&\text{otherwise}.
\end{cases}
\]
Mutation reverses the arrows incident with $k$. Checking the two possible orientations of each such edge, and using (R1) or (R3) for nonadjacent generators, gives, first on generators and hence on all elements,
\[
\begin{aligned}
\Phi_k\bigl(\widetilde\Phi_k(x)\bigr) & = s_kxs_k^{-1} &&(x\in M_\Gamma),\quad \widetilde\Phi_k\bigl(\Phi_k(y)\bigr) & = t_kyt_k^{-1} &&(y\in M_{\Gamma'}).
\end{aligned}
\]
Since $\Phi_k(t_k)=s_k$, we define
\[
\Psi_k:M_\Gamma\longrightarrow M_{\Gamma'}, \qquad \Psi_k(x)=t_k^{-1}\widetilde\Phi_k(x)t_k.
\]
For $x\in M_\Gamma$ and $y\in M_{\Gamma'}$, we have
\[
\begin{aligned}
\Phi_k\bigl(\Psi_k(x)\bigr) & =s_k^{-1}\Phi_k\bigl(\widetilde\Phi_k(x)\bigr)s_k =s_k^{-1}(s_kxs_k^{-1})s_k=x,\\
\Psi_k\bigl(\Phi_k(y)\bigr) &=t_k^{-1}\widetilde\Phi_k\bigl(\Phi_k(y)\bigr)t_k =t_k^{-1}(t_kyt_k^{-1})t_k=y.
\end{aligned}
\]
Thus $\Phi_k$ is an isomorphism with inverse $\Psi_k$.
\end{proof}

\section*{Acknowledgements}
This work was partially supported by the National Natural Science Foundation of China (Grant No. 12001254) and by Lanzhou University (No. 561120205).

\end{document}